\magnification=\magstep1

\def\double{\baselineskip=24truept plus0truept minus0truept}

\centerline{\phantom{STUFF}}

\bigskip

\centerline{\bf Normalization integrals of orthogonal Heun functions} 

\bigskip
\bigskip

\centerline{Peter A. Becker}

\bigskip
\bigskip

\centerline{\it Center for Earth Observing and Space Research}
\centerline{\it Institute for Computational Sciences and Informatics,}
\centerline{\it and Department of Physics and Astronomy,}
\centerline{\it George Mason University, Fairfax, VA 22030-4444}

\bigskip
\bigskip

\centerline{\tt (pbecker@gmu.edu)}

\bigskip
\bigskip

\double

A formula for evaluating the quadratic normalization
integrals of orthogonal Heun functions over the real interval
$0 \le x \le 1$ is derived using a simple limiting procedure
based upon the associated differential equation. The resulting
expression gives the value of the normalization integral explicitly
in terms of the local power-series solutions about $x=0$ and $x=1$
and their derivatives. This provides an extremely efficient
alternative to numerical integration for the development of
an orthonormal basis using Heun functions, because all of the
required information is available as a by-product of the search
for the eigenvalues of the differential equation.

\bigskip


\bigskip

Byline: Normalization integrals of Heun functions

\bigskip

({\it Accepted for publication in the
Journal of Mathematical Physics})

\vfil
\eject

\line{\bf I. INTRODUCTION \hfil}
\medskip

Heun's equation$^1$ is the most general Fuchsian equation
of second order with four regular singular points, and it
is therefore of considerable importance in mathematical physics.
Special
cases of the Heun equation include the
hypergeometric, confluent hypergeometric, Lam\'e, Bessel,
Legendre, and Laguerre equations.
As a practical matter, the most important solutions to the
Heun equation are those
orthogonal functions satisfying prescribed boundary conditions
at two adjacent singular points.$^2$ The development of an
orthonormal basis using these functions is a two-step process.
The first step is to search for the eigenvalues, and the second is
to evaluate the quadratic normalization integrals of the associated
eigenfunctions, which are orthogonal Heun functions. The normalization
integrals are usually evaluated numerically, which is not
an especially efficient procedure given the nature of the
eigenfunctions and the associated series representations.
In order to provide a useful alternative to numerical integration,
in this paper we derive
a new formula for directly evaluating the quadratic normalization
integrals of orthogonal Heun functions over the real interval $[0,1]$.
The formula obtained (eq. 29)
utilizes only information available as a by-product of the
search for the eigenvalues, and therefore greatly improves
both the efficiency and the accuracy of numerical procedures
involving Heun functions.

\medskip
\line{\bf II. HEUN FUNCTIONS \hfil}
\medskip

We begin by writing Heun's equation in the standard form first
adopted by Erd\'elyi~et al.,$^3$
$$
{d^2 y \over dx^2} + \left({\gamma\over x} + {\delta\over x-1}
+ {\epsilon\over x-a}\right) {dy \over dx} + {\alpha\beta x - \lambda
\over x(x-1)(x-a)}\,y = 0\,, \eqno(1)
$$
represented by the Riemann $P$-symbol
$$
P \left\{ \matrix{0&1&a&\infty\cr
                  0&0&0&\alpha&x\cr
                  1-\gamma&1-\delta&1-\epsilon&\beta\cr}\right\}\,,
\eqno(2)
$$
with regular singular points located at $x=0,1,a,\infty$.
The parameter $\lambda$ falls outside the domain of the
usual Riemann classification scheme, and is therefore referred to as
an accessory (or auxiliary) parameter. In many applications,
$\lambda$ plays the role of an eigenparameter. The five exponent
parameters $\alpha$, $\beta$, $\gamma$, $\delta$, $\epsilon$
are connected via Riemann's relation
$$
\alpha + \beta - \gamma - \delta - \epsilon + 1 = 0 \,, \eqno(3)
$$
and therefore only four of them are independent. The total number
of free parameters is six, and this number cannot be reduced by any
transformation.

Heun$^1$ used the method of Frobenius to derive local power-series
solutions to (1), generating a three-term recursion relation for
the expansion coefficients. Two linearly independent power-series
(Frobenius) solutions exist in the neighborhood of any one of the
singular points, and, in general, analytic continuation of a single
Frobenius solution about one singularity into the
neighborhood of a second, adjacent singularity generates
a linear combination of the two local Frobenius solutions about
the second singularity.

The most important solutions are those that are simultaneously
local Frobenius solutions about {\it two} adjacent singular points.
These are referred to as {\it Heun functions}, and often arise when
physical boundary conditions are applied to solutions of the
differential equation. Taking the parameters $\alpha$, $\beta$,
$\gamma$, $\delta$, $a$ to be constants, the problem of finding
a Heun function becomes a singular Sturm-Liouville eigenvalue
problem for $\lambda$. However, since no formula is available
for performing the analytic continuation between two adjacent
singularities in the case of the general Heun equation (in
contrast to the subcase of the hypergeometric equation), no
general closed-form expression for the eigenvalues $\lambda_n$
exists. Heun functions are infinite series in general, although
in special cases the series truncates, leaving a Heun polynomial.

There are four classes of Heun functions for a given pair of
adjacent singularities $x=s_0$ and $x=s_1$ that are distinguished
by the values of the corresponding exponents $(\sigma_0,\sigma_1)$.
In this paper, we set $s_0=0$ and $s_1=1$, and classify the Heun
functions according to the usual scheme based upon the values of
the associated exponents:
class I $(0,0)$; class II $(1-\gamma,0)$;
class III $(0,1-\delta)$; and class IV $(1-\gamma,1-\delta)$.
We focus here on the behavior of Heun functions
within the real interval $0 \le x \le 1$, and we assume throughout that
$a \not\in [0,1]$.

\medskip
\line{\bf III. DETERMINATION OF THE EIGENVALUES \hfil}
\medskip

Let the function $y_0(\lambda,x)$ be a local Frobenius solution
of (1) in the neighborhood of $x=0$, and let the function
$y_1(\lambda,x)$ be a local Frobenius solution in the neighborhood
of $x=1$. The properties of these solutions, including the recursion
relation for the expansion coefficients, have been fully discussed.$^2$
Let us suppose that these functions are normalized according to the
prescription adopted by Heun, so that
$$
\lim_{x \to 0} \  x^{-\sigma_0}\,y_0(\lambda,x) =
\lim_{x \to 1} \ (x-1)^{-\sigma_1}\,y_1(\lambda,x) = 1 \,.
\eqno(4)
$$
It is known from the basic theory of these solutions that
in general $y_0$ converges for $|x| < \min(1,|a|)$ and $y_1$
converges for $|x-1| < \min(1,|a-1|)$. Hence when $a \not\in [0,1]$
as assumed here, there exists a region
of {\it mutual convergence}, within which both $y_0$ and $y_1$
converge. The region of mutual convergence is contained within
the interval $[0,1]$.

The Wronskian of $y_0$ and $y_1$ is defined by
$$
W(\lambda,x) \equiv  y_0(\lambda,x){\partial y_1\over\partial x}
(\lambda,x) - y_1(\lambda,x){\partial y_0\over\partial x}
(\lambda,x)\,. \eqno(5)
$$
In order to obtain a Heun function, the Wronskian must vanish,
and therefore the eigenvalue equation for $\lambda_n$ becomes
$$
W(\lambda_n,x) = 0 \,. \eqno(6)
$$
When this condition is satisfied, $y_0(\lambda_n,x)$ and
$y_1(\lambda_n,x)$ are linearly dependent functions (although
not equal in general), and the solution to (1) with
$\lambda=\lambda_n$ is the Heun function $H_n(x)$. We set
the normalization of $H_n(x)$ by stipulating that
$H_n(x)=y_0(\lambda_n,x)$ in the neighborhood of $x=0$.
In the neighborhood of $x=1$, we find that
$H_n(x)=A(\lambda_n)\,y_1(\lambda_n,x)$, where the value
of $A(\lambda_n)$ is determined by requiring that
$H_n(x)$ be continuous at an arbitrary point within
the region of mutual convergence of $y_0$ and $y_1$.

We can establish the functional form of the Wronskian by examining
the self-adjoint version of the Heun equation,
$$
{\cal L}\,y - \lambda\,\omega(x)\,y = 0\,, \eqno(7)
$$
where the weight function $\omega(x)$ is defined by
$$
\omega(x) \equiv  x^{\gamma-1} (x-1)^{\delta-1}
(x-a)^{\epsilon-1}\,, \eqno(8)
$$
and the operator ${\cal L}$ is defined by
$$
{\cal L}\,y \equiv {d \over dx}\left[x^\gamma (x-1)^\delta
(x-a)^\epsilon \, {dy \over dx}\right] + \alpha\,\beta\, x \,
\omega(x)\,y \,.
\eqno(9)
$$
Since $y_0$ and $y_1$ are each solutions of (1) for the same value
of $\lambda$, we may write
$$
y_0\,\left[{\cal L} - \lambda\,\omega(x)\right]\,y_1
- y_1\,\left[{\cal L} - \lambda\,\omega(x)\right]\,y_0 = 0\,.
\eqno(10)
$$
This yields an equation for the Wronskian;
$$
-{1\over W}{d W\over dx} = {\gamma\over x} + {\delta\over x-1}
+ {\epsilon\over x-a} \,,
$$
with solution
$$
W(\lambda,x) = D(\lambda) \, x^{-\gamma} (x-1)^{-\delta}
(x-a)^{-\epsilon} \,, \eqno(11)
$$
where $D(\lambda)$ is an unknown function. Note that the
Wronskian vanishes for $D(\lambda_n) = 0$, which is
independent of $x$. Hence when we evaluate $W(\lambda,x)$
using (5) in order to calculate the eigenvalues using (6),
we are free to pick any convenient value for $x$ that lies
within the region of mutual convergence of $y_0$ and $y_1$.

\medskip
\line{\bf IV. ORTHOGONALITY RELATIONS \hfil}
\medskip

Let $H_n(x)$ and $H_m(x)$ be Heun functions of the same class
associated with
eigenvalues $\lambda_n$ and $\lambda_m$. Since they are each
solutions of (1), we have, by analogy with~(10),
$$
H_n\,\left[{\cal L} - \lambda_m\,\omega(x)\right]\,H_m
- H_m\,\left[{\cal L} - \lambda_n\,\omega(x)\right]\,H_n = 0\,.
\eqno(12)
$$
After simplifying and integrating over the interval $[0,1]$,
we obtain
$$
(\lambda_n - \lambda_m)\,\int_0^1\omega(x)\,H_n(x)\,H_m(x)\,dx
= p(x)\,
\left(H_m\,{dH_n\over dx} - H_n\,{dH_m\over dx}
\right)\bigg\vert_0^1\,, \eqno(13)
$$
where
$$
p(x) \equiv x^\gamma (x-1)^\delta (x-a)^\epsilon \,. \eqno(14)
$$
Based upon the asymptotic behavior of the local Frobenius solutions
$y_0$ and $y_1$, we find that the right-hand side of (13) vanishes
when one of the following sets of class-dependent conditions is
satisfied;
$$
\matrix{$class I:$~~~   &\Re\gamma > 0\,,   &\Re\delta > 0\cr
        $class II:$~~   &\Re\gamma < 2\,,   &\Re\delta > 0\cr
        $class III:$\,  &\Re\gamma > 0\,,   &\Re\delta < 2\cr
        $class IV:$\,   &\Re\gamma < 2\,,   &~\,\Re\delta < 2\,.\cr} \eqno(15)
$$
We shall refer to (15) as the set of {\it existence
conditions} for orthogonal Heun functions on the interval $[0,1]$.
When these conditions are met, we obtain the standard
orthogonality relation
$$
(\lambda_n - \lambda_m)\,\int_0^1\omega(x)\,H_n(x)\,H_m(x)\,dx
= 0\,. \eqno(16)
$$

\medskip
\line{\bf V. NORMALIZATION INTEGRALS \hfil}
\medskip

Naturally, the integral in (16) does not vanish when $n = m$,
and in this case it is necessary to establish its value in order
to develop an orthonormal basis using orthogonal Heun functions.
Numerical integration is always available as an option,
but this is a very inefficient approach to the problem. Other
procedures have been devised, the most interesting being the
method developed by Erd\'elyi,$^4$ which is based upon expansions
of Heun functions as series of degenerate hypergeometric functions
(Jacobi polynomials). Lambe and Ward$^5$ developed a technique for
evaluating normalization integrals for Heun polynomials, but these
results are not applicable to the more general (and much more common)
case of Heun functions. In this section we develop a new formula
for the explicit evaluation of these integrals.

We can derive a formula for evaluating the quadratic normalization
integral
$$
I_n \equiv \int_0^1 \omega(x)\,\left[H_n(x)\right]^2\,dx
\eqno(17)
$$
by generalizing
the approach taken in the preceding section.
Proceeding as before, we note that since $y_0(\lambda,x)$ and
$y_0(\lambda_n,x)$ are each Frobenius solutions of (1) in the
neighborhood of $x=0$, it follows by analogy with (12) that
$$
y_0(\lambda,x)\,[{\cal L}-\lambda_n\,\omega(x)]\,y_0(\lambda_n,x)
- y_0(\lambda_n,x)\,[{\cal L}-\lambda\,\omega(x)]\,y_0(\lambda,x)
= 0\,, \eqno(18)
$$
where $\lambda_n$ is an eigenvalue and $\lambda$ is arbitrary.
After simplifying and integrating with respect to $x$, we now obtain
$$
\eqalign{
(\lambda - \lambda_n)\,\int_0^x
&\omega(x')\,y_0(\lambda,x')\,
y_0(\lambda_n,x')\,dx'\cr
& = p(x)
\left[y_0(\lambda_n,x)\,{\partial y_0\over \partial x}(\lambda,x)
- y_0(\lambda,x)\,{\partial y_0\over \partial x}(\lambda_n,x)
\right]
\,,\cr} \eqno(19)
$$
where the right-hand side vanishes as $x \to 0$ provided the
appropriate set of existence conditions in (15) is satisfied.

Upon examination of (19), we find that both sides of the equation
vanish as $\lambda \to \lambda_n$. We can therefore establish the
value of the indefinite integral in the limit $\lambda \to \lambda_n$
using L'H\^opital's rule, which yields
$$
\eqalign{
\int_0^x
\omega(x')\,[&y_0(\lambda_n,x')]^2\,dx'\cr
& = p(x)
\left[y_0(\lambda_n,x)\,{\partial^2 y_0\over \partial \lambda
\partial x}(\lambda_n,x)
- {\partial y_0 \over \partial \lambda}(\lambda_n,x)\,
{\partial y_0\over \partial x}(\lambda_n,x)\right]
\,.\cr} \eqno(20)
$$
To obtain a formula for the desired normalization integral
$I_n$, we must let $x \to 1$ in (20), and this requires analytic
continuation of $y_0$ into the neighborhood of $x=1$.
For general values of $\lambda$, we write the analytic continuation
of $y_0$ as
$$
y_0(\lambda,x) = A(\lambda)\,y_1(\lambda,x)
+ B(\lambda)\,\tilde y_1(\lambda,x)\,, \eqno(21)
$$
where $y_1$ is the Frobenius solution about $x=1$ with exponent zero,
and $\tilde y_1$ is the Frobenius solution about $x=1$ with exponent
$1-\delta$. At this point we shall restrict our attention to Heun
functions of classes I or II, so that the exponent at $x=1$ is
zero. This restriction will be removed later.
Note that $B(\lambda_n)$ must vanish so that we obtain
a class I or II Heun function when $\lambda=\lambda_n$. Substituting
into (20), we obtain after some algebra
$$
\eqalign{
I_n &\equiv
\int_0^1 \omega(x)\,[H_n(x)]^2\,dx\cr
& = \lim_{x \to 1} \left[
A\,y_1
\left(A\,{\partial^2 y_1 \over \partial\lambda\partial x}
+ {dB\over d\lambda}\,{\partial\tilde y_1\over\partial x}\right)
- A\,{\partial y_1 \over\partial x}
\left(A\,{\partial y_1 \over \partial\lambda}
+ {dB\over d\lambda}\,\tilde y_1\right)
\right] p \ 
\bigg\vert_{\lambda=\lambda_n}\,,\cr} \eqno(22)
$$
where we have used the fact that in the neighborhood of $x=1$,
$H_n(x)$ is the analytic continuation of $y_0(\lambda_n,x)$, and
in the neighborhood of $x=0$, $H_n(x)$ is identical to
$y_0(\lambda_n,x)$.

Based upon asymptotic analysis of $y_1$ and $\tilde y_1$, we find
that in the limit $x \to 1$, all of the terms on the right-hand
side of (22) vanish except the second one, so that we are left with
$$
I_n = \lim_{x \to 1}\ \ 
p(x)\,A(\lambda_n)\,{dB\over d\lambda}(\lambda_n)\,
{\partial\tilde y_1\over\partial x}(\lambda_n,x)\,
y_1(\lambda_n,x)\,.\eqno(23)
$$
Hence we need only evaluate ${dB\over d\lambda}(\lambda_n)$ and
$A(\lambda_n)$ in terms of known functions in order to obtain our
final result for the normalization integral $I_n$.
Evaluation of $A(\lambda_n)$ is a simple matter, since
the continuity of $H_n(x)$ requires that
$$
A(\lambda_n) = {y_0(\lambda_n,x) \over y_1(\lambda_n,x)}
\eqno(24)
$$
for any $x$ within the region of mutual convergence of $y_0$ and $y_1$.
Evaluation of ${dB\over d\lambda}(\lambda_n)$ is slightly more
complicated. We begin by differentiating (21) with respect to
$x$ to obtain
$$
{\partial y_0 \over \partial x}(\lambda,x)
= A(\lambda)\,{\partial y_1 \over \partial x}(\lambda,x)
+ B(\lambda)\,{\partial \tilde y_1 \over \partial x}(\lambda,x)\,.
\eqno(25)
$$
Solving (21) and (25) for $B(\lambda)$ yields
$$
B(\lambda) = {
y_0\,{\partial y_1\over\partial x}
- y_1\,{\partial y_0\over\partial x}
\over
\tilde y_1\,{\partial y_1\over\partial x}
- y_1\,{\partial\tilde y_1\over\partial x}
}\,\ \ . \eqno(26)
$$
Bearing in mind that $B(\lambda_n)=0$, we find upon differentiating
(26) with respect to $\lambda$ that
$$
{dB\over d\lambda}(\lambda_n) = {
{\partial W\over\partial\lambda}
\over
\tilde y_1\,{\partial y_1\over\partial x}
- y_1\,{\partial\tilde y_1\over\partial x}
}\,\,\bigg\vert_{\lambda=\lambda_n}\,, \eqno(27)
$$
where the Wronskian $W(\lambda,x)$ is defined by (5).
We can use (24) and (27) respectively to replace $A(\lambda_n)$
and ${dB\over d\lambda}(\lambda_n)$ in (23), yielding
$$
I_n = \lim_{x \to 1}\ \ p\,{y_0\over y_1}\,{
{\partial W\over\partial\lambda}
\over
\tilde y_1\,{\partial y_1\over\partial x}
- y_1\,{\partial\tilde y_1\over\partial x}
}\,
{\partial\tilde y_1\over\partial x}\,y_1\,\,
\bigg\vert_{\lambda=\lambda_n}\,.
\eqno(28)
$$
The most singular term in the denominator is the one containing
${\partial\tilde y_1\over\partial x}$, and therefore in the limit
$x \to 1$ our final result for the normalization integral becomes
$$
\int_0^1 \omega(x)\,\left[H_n(x)\right]^2\,dx
= - p(x)\,{\partial W\over\partial\lambda}(\lambda_n,x)\,
{y_0(\lambda_n,x) \over y_1(\lambda_n,x)}\,. \eqno(29)
$$
In passing to this expression, we have made use of
the fact that ${\partial W\over\partial\lambda}(\lambda_n,x)\,p(x)$
and $y_0(\lambda_n,x)/y_1(\lambda_n,x)$ are both independent of
$x$. Hence the right-hand side of (29) is actually an {\it invariant},
which may be evaluated at any point within the region of mutual
convergence of the local power-series solutions $y_0$
and $y_1$.

\medskip
\line{\bf VI. EXTENSION TO HEUN FUNCTIONS OF CLASSES III AND IV \hfil}
\medskip

We have derived our main result (29) under the assumption
that $H_n(x)$ is a class I or II Heun function. In the case of
a Heun function of class III or IV, the exponent at $x=1$ is
$1-\delta$, and the derivation presented above must be modified
slightly. First we note that (20) remains valid because it is
written in terms of the local solution $y_0$ about $x=0$. The
analytic continuation of $y_0$ into the neighborhood of $x=1$
is still performed using (21), except now we interchange
the definitions of $y_1$ and $\tilde y_1$, so that $y_1$ is the
local Frobenius solution about $x=1$ with exponent $1-\delta$
and $\tilde y_1$ is the local solution about $x=1$ with exponent
zero.
In this case we again require $B(\lambda_n)=0$, and following
the same procedure as before we regain expression (22). Due to
the fact that the definitions of $y_1$ and $\tilde y_1$ have
been interchanged, in the limit $x \to 1$, we find by asymptotic
analysis that only the {\it final} term in (22) now contributes,
so that in this case we obtain
$$
I_n = \lim_{x \to 1}\ \ 
- p(x)\,A(\lambda_n)\,{dB\over d\lambda}(\lambda_n)\,
{\partial y_1\over\partial x}(\lambda_n,x)\,
\tilde y_1(\lambda_n,x)\,.\eqno(30)
$$
Evaluation of $A(\lambda_n)$ and ${d B\over d\lambda}(\lambda_n)$
proceeds exactly as before, and we simply regain expressions
(24) and (27). Using these results in (30) now yields
$$
I_n = \lim_{x \to 1}\ \ - p\,{y_0\over y_1}\,{
{\partial W\over\partial\lambda}
\over
\tilde y_1\,{\partial y_1\over\partial x}
- y_1\,{\partial\tilde y_1\over\partial x}
}\,
{\partial y_1\over\partial x}\,\tilde y_1\,\,
\bigg\vert_{\lambda=\lambda_n}\,.
\eqno(31)
$$
This is similar to (28), except now the most singular term in the
denominator is the one containing ${\partial y_1\over \partial x}$,
and therefore in the limit $x \to 1$ we obtain
$$
I_n
= - p(x)\,{\partial W\over\partial\lambda}(\lambda_n,x)\,
{y_0(\lambda_n,x) \over y_1(\lambda_n,x)}\,,
$$
which is identical to (29). We are therefore led to the following
general conclusion: {\it Equation (29) holds for Heun functions
of any class, provided the existence conditions (15) are satisfied}.
We discuss the significance of this result and its natural role
in computational algorithms below.

\medskip
\line{\bf VII. COMPUTATIONAL IMPLICATIONS \hfil}
\medskip

Equation (29) has important
practical consequences for the design of numerical algorithms used
to develop orthonormal systems based on Heun functions, which are
the solutions of interest in many mathematical and physical
situations.
In such cases the first step towards solution
is the search for the associated eigenvalues $\lambda_n$. This search
must proceed numerically in general, due to the lack of a closed-form
expression for $\lambda_n$ in terms of the parameters $\alpha$,
$\beta$, $\gamma$, $\delta$, $a$. The eigenvalues
are obtained by isolating the roots of the Wronskian in (6).
The root finding usually proceeds via Newton's method
or possibly some more sophisticated algorithm. Most of these
techniques require the evaluation of $W$ and
${\partial W\over\partial\lambda}$ in order to generate a
revised estimate of the true root $\lambda_n$, and
the evaluation of these functions in
turn involves the determination of the quantities
$$
y_0\ \ , \ \
{\partial y_0\over\partial x}\ \ , \ \
{\partial y_0\over\partial\lambda}\ \ , \ \
{\partial^2 y_0\over \partial \lambda\partial x}\ \ , \ \
y_1\ \ , \ \
{\partial y_1\over\partial x}\ \ , \ \
{\partial y_1\over\partial\lambda}\ \ , \ \
{\partial^2 y_1\over \partial \lambda\partial x}\ \ , \ \
$$
at each iteration.
The values of $y_0$ and $y_1$ can be obtained using the well-known
power-series representations, and the values of the derivatives can
be obtained using term-by-term differentiation. It is
straightforward to demonstrate that the radii of convergence of the
series for the derivatives are identical to those for the corresponding
fundamental series, which are discussed in \S~II.

Once the eigenvalues $\lambda_n$ have been determined to acceptable
precision,
one generally needs to evaluate the associated quadratic normalization
integrals $I_n$ in order to develop a set of orthonormal
basis functions $h_n(x)$
using
$$
h_n(x) \equiv {H_n(x) \over I_n^{1/2}} \,, \eqno(32)
$$
with normalization
$$
\int_0^1 \omega(x)\,[h_n(x)]^2\,dx = 1 \,. \eqno(33)
$$
The conventional approach to the problem of evaluating $I_n$
is to integrate (17) numerically, in which case thousands of
evaluations of $H_n(x)$ would generally be required in order
to establish the value of $I_n$ to reasonable accuracy.
However, such an inefficient procedure is no longer necessary
with the availability of (29), because it allows the determination
of $I_n$ to high precision using only the values of
$y_0$, $y_1$, and ${\partial W\over\partial\lambda}$ obtained
in the final iteration of the root-finding stage of the algorithm.
Hence no substantial additional calculation is necessary in order
to determine $I_n$.
The computational time required for
developing orthonormal systems of Heun functions can therefore
be reduced by several orders of magnitude by using (29) instead
of numerical integration.

We close by making a comparison between the method for evaluating
$I_n$ outlined here and that suggested by Erd\'elyi,$^4$
which utilizes Svartholm's$^6$ expansions of Heun functions as
series of degenerate hypergeometric functions, essentially Jacobi
polynomials. First of all, it is worth noting that in Erd\'elyi's
method, the root-finding approach is the same as that outlined
above, because his procedure assumes prior knowledge of the
eigenvalues. With the eigenvalues already determined, the
calculation of the coefficients for the Jacobi expansion
proceeds via a three-term recursion relation similar to that
derived by Heun for the coefficients of the power-series expansion.
Using the familiar result for the quadratic normalization
integrals of the Jacobi polynomials, it is a simple matter
to determine $I_n$ from the coefficients of the Jacobi
expansion. While Erd\'elyi's method is interesting from the point
of view of functional analysis, it is obvious that his procedure
entails much more computation than the evaluation of $I_n$
using (29), which we again emphasize requires nothing more than
information available as a by-product of the search for the
eigenvalues. In conclusion, we point out that the application
of L'H\^opital's rule used here to evaluate the normalization
integrals of Heun functions can also be used to obtain similar
results for more general Sturm-Liouville problems.

\medskip
\centerline{\bf REFERENCES}
\medskip

$^1$K. Heun, {\it Zer Theorie der Riemann'schen Functionen
zweiter Ordnung mit vier Verzweigungspunkten}, Math. Ann., 33,
pp. 161--179 (1889).

$^2$A. Ronveaux, ed., {\it Heun's Differential Equations}
(Oxford University Press, New York, NY, 1995).

$^3$A. Erd\'elyi, W. Magnus, F. Oberhettinger, and F. G. Tricomi,
{\it Higher Transcendental Functions}, vol. III (McGraw-Hill,
New York, NY, 1955).

$^4$A. Erd\'elyi, {\it Certain expansions of solutions of the
Heun equation}, Quart. J. Math. Oxford Ser., 15, pp. 62--69 (1944).

$^5$C. G. Lambe and D. R. Ward, {\it Some differential equations
and associated integral equations}, Quart. J. Math. Oxford Ser.,
5, pp. 81--97 (1934).

$^6$N. Svartholm, {\it Die L\"osung der Fuchsschen
Differentialgleichung zweiter Ordnung durch hypergeometrische
Polynome}, Math. Ann., 116, pp. 413--421 (1939).

\bye